\documentclass{cseg2015}
\usepackage{xy}
\usepackage{amsmath,graphicx,natbib}
\usepackage{rotating}
\usepackage{lscape}
\definecolor{burntorange}{RGB}{225,100,0}
\usepackage{setspace}
\graphicspath{ {figs/}}

\begin{document}
 \vspace*{-2.cm}
\begin{center}
\includegraphics[height=3.2cm]{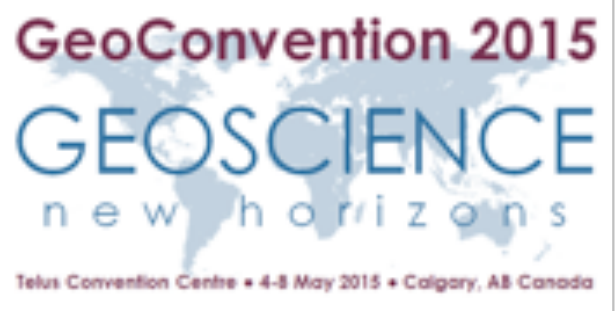}
\end{center}

{\LARGE {\bf Block row recursive least squares migration}} 

 \vspace*{-0.1cm}

{\it Nasser Kazemi and Mauricio Sacchi}

\vspace*{-0.2cm}

{\it University of Alberta, Canada}

\section{Summary}
Recursive estimates of large systems of equations in the context of least squares  fitting is a common practice in different fields of study. For example, recursive adaptive filtering is extensively used in signal processing and control applications. The necessity of solving least squares problem via recursive algorithms comes from the need of fast real-time signal processing strategies. Computational cost of using least squares algorithm could also limits the applicability of this technique in geophysical problems. In this paper, we consider recursive least squares solution for wave equation least squares migration with sliding windows involving several rank K downdating and updating computations. This technique can be applied for dynamic and stationary processes. One can show that in the case of stationary processes, the spectrum of the preconditioned system is clustered around one and the method will converge superlinearly with probability one, if we use enough data in each windowed setup. Numerical experiments are reported in order to illustrate the efectiveness of the technique for least squares migration.

\section{Introduction}
 Recursive least squares algorithms are extensively used for dynamic system of equations (e.g., radar data) and also for adaptive filtering in the case of dynamic and stationary environments. The main idea is to implement infinite memory algorithms by solving the problem via introducing one data point at a time to the system of equations. \cite {Mostafa} used rank one update of the recursive least squares fitting with some exponentially weighted forgetting factor in the context of $f-x$ adaptive filtering for seismic interpolation. However, very little effort has been made in the geophysical community to adapt recursive least squares techniques for different practical applications. This is mainly because, recursive least squares algorithms require explicit matrices. In large scale problems (e.g., migration) every thing will be done on the fly and also adding or removing one row from the system of equations (e.g., de-migration operator) has no physical meaning. Moreover, the one data point update scheme in recursive least squares approach causes stability issues. In other words, it is difficult to track the stability of the process. To tackle these short comings, one can update the system of equations in a block fashion by introducing more than one data point at each step, making the algorithm faster and more stable. In this paper we will try to adapt these kind of techniques for least squares migration.  
 
\section {Least squares Migration}

Data generating model under the action of de-migration operator ${\bf A}$ can be written as
\begin{equation}\label{eq:1}
{\bf d}={\bf A} {\bf m}+{\bf n},
\vspace{-.2 cm}
\end{equation}
where ${\bf d}$ is $N\times 1$ vectorized version of the recorded data at the surface, ${\bf m}$ is the migrated model with $M\times 1$ and ${\bf n}$ is the noise content and, often, also a term to absorb waves not modelled by
the de-migration operator. Using the adjoint operator of de-migration, one can estimate the migrated model
\begin{equation}\label{eq:2}
\hat {\bf m}={\bf A}^T {\bf d},
\vspace{-.2 cm}
\end{equation}
where $\hat {\bf m}$ is the adjoint estimated image of the earth and ${\bf A}^T$ is the migration operator.
While $\hat{\bf{m}}$ can capture the main structures of the true geological model $\bf{m}$, the produced model does not honour the data. In other words, application of the de-migration operator on the migrated image yields poor data prediction.  In addition, the migrated image contains blurring and sampling artifacts. These artifacts come from the fact that migration and de-migration operators are not orthogonal and the energy of the signal in the complimentary image space of the operator will be zeroed out. To tackle the problem, we will solve
\begin{equation}\label{eq:3}
  \tilde{\bf{m}} = \underset{{\bf m}}{\operatorname{argmin} }\; ||{\bf{A}}\;{\bf{m}}-{\bf{d}}||_2^2+\lambda\;||{\bf{m}}||_2^2,
  \vspace{-.2 cm}
\end{equation}
where ${\bf m}$ is desired model, ${\bf{A}}$ is de-migration operator,  ${\bf{d}}$ is recorded data and $\lambda$ is a regularization parameter \cite[]{Nemeth,  Kuhl, Kaplan2}.  The cost function of equation $\ref{eq:3}$ is convex and has a close form solution
\begin{equation}\label{eq:4}
\tilde{\bf{m}}=({\bf{A}}^T{\bf{A}}+\lambda\;{\bf{I}})^{-1}\;{\bf{A}}^T{\bf{d}}.
\vspace{-.2 cm}
\end{equation}
However, the computational cost of least squares approach is very high and to make the algorithm faster, one can approximate the Hessian inverse \cite[]{Hu,kazemi}.
 In next section we will propose a sliding window scheme for solving equation \ref{eq:3}  in recursive fashion. The main motivations behind the method are: reducing the computational cost and at the same time producing migrated images that honour the recorded data using memory limited resources.
 
 \section{Block row recursive least squares migration}
\vspace{-.3 cm}
 In this section we will follow the recursive least squares solution via rank $K$ updating and rank $K$ downdating procedure introduced by \cite{Ng}. However, there are some differences between the proposed method in \cite{Ng} with our technique. We are not considering near Toeplitz structure for data matrix and also in the updating procedure, we just use previous solution of the block row setup as an initial solution for the next sliding window. This way we will give up fast implementation of the technique in the favour of not considering special structures in the de-migration operator. 

To explain the block row recursive least squares method, let us consider again the problem of equation \ref{eq:3}. In recursive least squares computations, it is required to calculate ${\bf m}$ while observations are successively added to, or deleted from the system of equations. 
Suppose we have estimated the model with first set of measurement points ${\bf d}_0$ in least squares sense
\begin{equation}\label{eq:5}
{\bf m}_0 = ({\bf{A}}_0^T{\bf{A}}_0)^{-1}\;{\bf{A}}_0^T{\bf{d}}_0,
\end{equation}
Now, the question we should ask is that by introducing new data points to the system of equations can the best estimate for the combined system ${\bf A}_0\;{\bf m}={\bf d}_0$ and ${\bf A}_1\;{\bf m}={\bf d}_1$ be estimated using only ${\bf m}_0$ and ${\bf d}_1$?

To do that, we define new matrix
\begin{equation}\label{eq:6}
{\bf P}_1^{-1} = {\bf{A}}_0^T{\bf{A}}_0+{\bf{A}}_1^T{\bf{A}}_1,
\vspace{-.1 cm}
\end{equation}
then we have 
\begin{equation}\label{eq:7}
{\bf m}_1 = {\bf P}_1\;({\bf{A}}_0^T{\bf{d}}_0+{\bf{A}}_1^T{\bf{d}}_1),
\vspace{-.1 cm}
\end{equation}
note that ${\bf m}_1$ is the best model for combined system of equations. At this point we need to eliminate ${\bf d}_0$ term from equation \ref{eq:7}. Let us rewrite equation \ref{eq:6} as
\begin{equation}\label{eq:8}
{\bf P}_1^{-1} = {\bf P}_0^{-1}+{\bf{A}}_1^T{\bf{A}}_1,
\vspace{-.1 cm}
\end{equation}
and after few algebraic manipulations, one can show that equation \ref{eq:7} changes to
\begin{equation}\label{eq:9}
{\bf m}_1 = {\bf m}_0+{\bf P}_1{\bf{A}}_1^T({\bf{d}}_1-{\bf{A}}_1{\bf{m}}_0),
\vspace{-.1 cm}
\end{equation}
where ${\bf P}_1{\bf{A}}_1^T$ is gaining factor. In the case of one data point update we can use Matrix Inversion Lemma (MIL) and calculate the gaining factor without direct inversion and in the case of block wise update (More than one data point update) to apply fast calculations of gaining matrix we need to explore MIL with QR decomposition. However, none of these approaches are applicable to least squares migration. It is mainly because, in migration every thing will be done on the fly and there are no simple explicit matrices for migration. 

\begin{figure}
\centerline{\includegraphics[width=1.0\columnwidth]{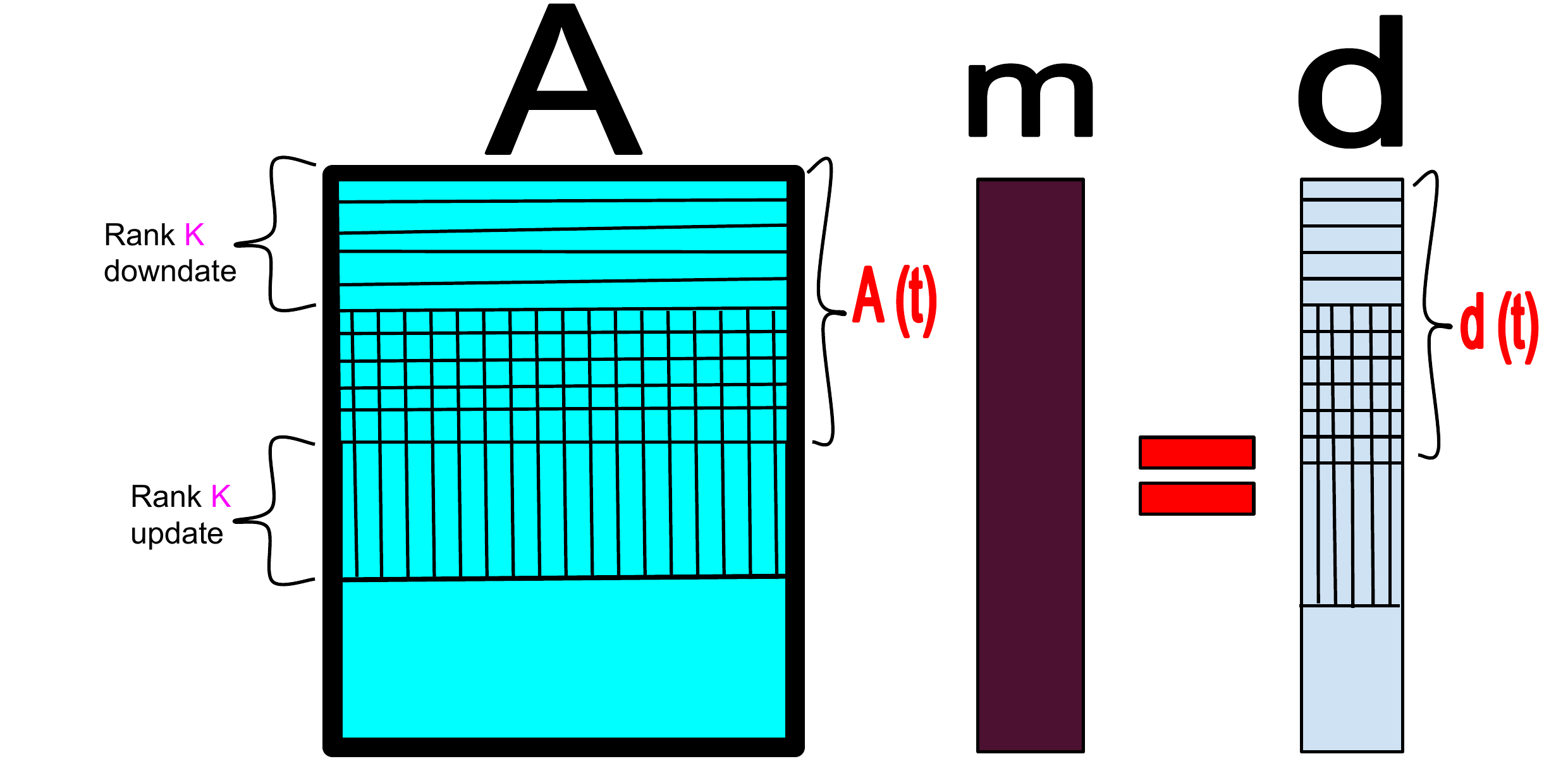}}

\caption{Schematic representation of the block row recursive least squares algorithm.   }
\vspace{-.1 cm}
\label{fig:1}
\end{figure}

Fortunately, \cite{Ng} showed that if the system of equations satisfy some assumptions,  we can relax the gaining factor term of equation \ref{eq:9} and solve the problem recursively on the overlapping windows using CG method.
Let us explain the method step by step. The least squares estimator at step $t$ can be found by solving for the $M\times 1$ vector ${\bf m}(t)$ in
\begin{equation}\label{eq:10}
 \tilde{\bf{m}}(t) = \underset{{\bf m}(t)}{\operatorname{argmin} }\; ||{\bf{A}}(t)\;{\bf{m}}(t)-{\bf{d}}(t)||_2^2+\lambda\;||{\bf{m}}(t)||_2^2,
\end{equation}
where $\tilde{\bf{m}}(t)=[\tilde{m}_1(t),\tilde{m}_2(t),...,\tilde{m}_M(t)]^T$ is the least squares solution of the model at step $t$ and ${\bf{d}}(t)=[d(t-Q+1),d(t-Q+2),...,d(t)]^T$ is the recorded data at step $t$ with size $Q\times 1$ and {\bf{A}}(t) is the $Q\times M$ data matrix at step $t$ and $Q$ is the length of sliding window. We consider causal systems (i.e., $t \geq Q$). Figure \ref{fig:1} shows the schematic representation of the setup. 
To update the solution recursively, we will add $K$ data points to the system of equations and remove $K$ data points from the beginnings of the previous data vector. We call this step as rank $K$ updating and downdating (see Figure \ref{fig:1} for more information). For this new configuration we have $t= Q+K$ and \cite{Ng} showed that one can use
 \begin{equation}\label{eq:11}
\tilde{\bf m}(K+Q) = \tilde{\bf m}(Q)+({\bf{A}}^T(K+Q){\bf{A}}(K+Q))^{-1}{\bf{A}}^T(K+Q)({\bf{d}}(K+Q)-{\bf{A}}(K+Q)\tilde{\bf{m}}(Q)),
\end{equation}
to estimate the best model in least squares sense that fits the whole system of equations at step $t=K+Q$. The second term in equation \ref{eq:11} can be interpreted as the least squares solution of the unpredicted part of the new data set by the solution of previous step. This warm start for the new sliding window will result in fast convergence of the CG method. This is mainly because the nearby data points are highly correlated in seismic acquisition. It is obvious that $t=Q,Q+K,Q+2K,...$ are step values at each rank $K$ updating and downdating setup. It is also worth mentioning that in each setup we must use proper pre-conditioners and regularization term to cluster the eigenvalues of the partial Hessians around 1. \cite {Ng} proved the convergence of this recursive least squares technique in the probabilistic terms. They showed that the method will converge superlinearly with probability 1 provided the underlying process satisfies some assumptions. First of all, the input discrete- time stochastic process should be stationary. Secondly, the auto-covariances of the kernels in each step should be absolutely summable. This in turn will assure the invertibility of the processes in each windowed setup. Thirdly, the variances between the auto-covariances of the kernels between different setups should be bounded. Finally, the stationary process has zero mean. All of these assumptions are valid for many time series analysis problems but it is not clear if they are fully applicable to migration problems. Nevertheless, the ideas in Ng and Plemmons (1996) can be explored for migration applications. In next section we will show the efficiency of the proposed method using Marmousi model.

\section {Examples}

In the following examples we used shot profile wave equation adjoint and forward operators with the split step correction in the context of adjoint, least squares and block row recursive least squares algorithms. 

\begin{figure}
\vspace{-1.5 cm}
\centerline{\includegraphics[width=.98\columnwidth]{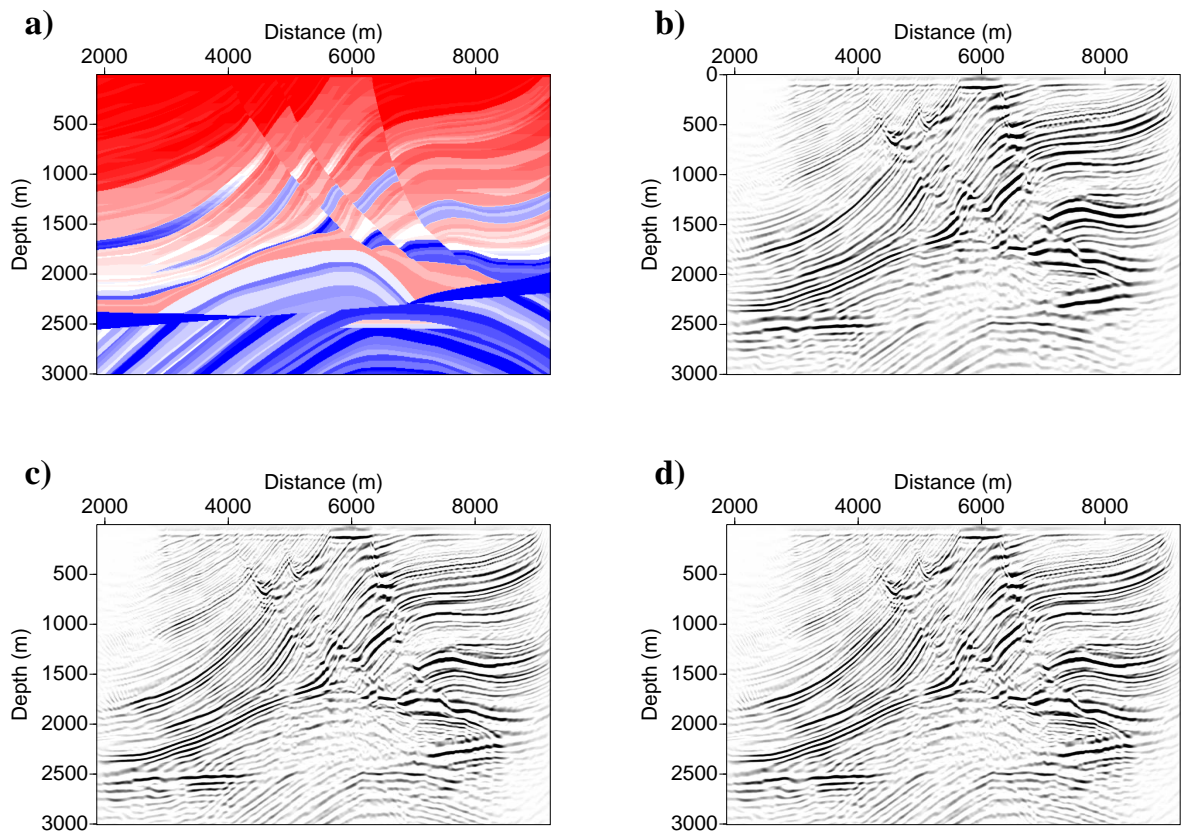}}

\caption{Comparison of the performances of different methods. a) True velocity field. b) Adjoint migrated model. c) Least squares migrated model. d) Block row recursive least squares migrated model.  }
\vspace{-.2 cm}
\label{fig:2}
\end{figure}

To test the performance of the proposed method, we applied the method on the Marmusi model. It is worth mentioning that the data set generated by finite difference method with a ricker wavelet with dominant frequency of $20\; Hz$. The data set consist of 240 shot gathers with $25m$ shot interval that modelled with an off end survey with receivers to the left of the source being pulled towards the right. Non-smooth velocity model of Marmousi is shown in Figure \ref{fig:2}a. Figure \ref{fig:2}b shows the adjoint migrated image of Marmousi data set using shot profile wave equation migration with split step correction. Figure \ref{fig:2}c shows the least squares migrated image of Marmousi data set after $8$ iterations using split step Fourier migration and de- migration operators as an adjoint and forward operators required by CG algorithm. We used some predefined threshold for the normal equation as a convergence criterion. Finally, Figure \ref{fig:2}d shows the migrated image produced by the block row recursive approach.  In the case of the block row recursive approach we tested different configurations and finally we used $5$ consecutive shot gathers in each group and we deleted $3$ shot gathers from the beginnings of the previous windowed setup and added $3$ new shot gathers to the end of new windowed setup. It is worth mentioning that we got the same results using different configurations. The block row recursive approach did a good job in recovering the true reflectivity model and the result is comparable to that of the least squares. Comparing the proposed 
method's result with adjoint, the block row recursive approach did a reasonable job in preserving the amplitude of the reflectors and removed some of the defocusing problems inherent in the adjoint migrated image. 

\section {Conclusion}

We present a block row recursive least squares migration method. This method uses block wise update of the de-migration operator via rank $K$ update and downdate in each set up while the new data points are successively added to the data vector. In each windowed set up, CG algorithm used to solve the system of equations in least squares sense. To have fast convergence, previous solution of the method implemented as an initial solution for the next block. This warm start will result in fast convergence of the CG algorithm. This is supported by the fact that nearby blocks are highly correlated. The results of applying this technique on Marmousi model convinced us that the block row recursive method can be a practical tool for improving the spatial resolution of the migrated images.

\section{Acknowledgment}
We thank the sponsors of the Signal Analysis and Imaging Group (SAIG) at the University of Alberta.
\bibliographystyle{seg}  
\bibliography{paper}
\end{document}